\documentclass[11pt, a4paper]{article}

\usepackage{amsmath, amssymb, amscd, amsfonts}

\usepackage{theorem}

\input xy
\xyoption{all}

\newtheorem{teo}{Theorem}[section]

\newtheorem{cor}[teo]{Corollary}

\newtheorem{prop}[teo]{Proposition}

\theorembodyfont{\normalfont}

\newtheorem{defi}[teo]{Definition}

\newtheorem{rmk}[teo]{Remark}

\newtheorem{es}[teo]{Example}

\newenvironment{proof}[1][]{\noindent\textbf{Proof#1}.  }{\begin{flushright} $\blacksquare$ \end{flushright}}

\numberwithin{equation}{section}

\newcommand {\fs}{\mathcal{O}}

\newcommand {\p}{\mathbb{P}^1}

\newcommand {\fr}{\rightarrow}

\renewcommand {\cal}{\mathcal}

\newcommand {\Sym}{\operatorname{Sym}} 

\newcommand {\sym}{\operatorname{Sym}}

\begin{document}

\title{Slope equalities for genus 5 surface fibrations}
\author{Elisa Tenni}
\date{}
\maketitle
\begin{abstract}\noindent
In \cite{AK} and \cite{K} Konno proved a \emph{slope equality} for fibred surfaces with fibres of odd genus and general fibre of maximal gonality. More precisely he found a relation between the invariants of the fibration and certain weights of special fibres (called the \emph{Horikawa numbers}). We give an alternative and more geometric proof in the case of a genus 5 fibration, under generality assumptions. In our setting we are able to prove that the fibre with positive Horikawa numbers are precisely the trigonal ones, we compute their weights explicitly and thus we exhibit explicit examples of regular surfaces with assigned invariants and Horikawa numbers.
\end{abstract}
\section{Introduction}
Let $f\!: S \fr B$ be a surface fibration. Suppose that it satisfies a prescribed condition (*) on the general fibre, and that there exists a rational number $\lambda$ and well-defined nonnegative rational numbers $\cal{H}(f,F)$ depending on the fibration $f$ and on the fibre $F$, and  vanishing on the fibres satisfying the general condition (*) such that the following relation is true for any fibration of the class:
\[K_f^2= \lambda \chi_f + \sum_{P \in B} \cal{H}(f, F_P)\]
Suppose moreover that there are infinitely many fibrations which satisfy the lower bound $K_f^2= \lambda \chi_f$.

In this situation, following the notation of Ashikaga, Konno and Reid in \cite{As}, \cite{K} and \cite{Reid2}, the above relation is called a \textit{slope equality}. The fibres with positive index $\cal{H}$ are called the \textit{atoms} of the fibration, and the numbers $\cal{H}(f,F)$ the \emph{Horikawa numbers} (or \emph{Horikawa indexes}) of the fibres.

The essence of a slope equality is that it is a relation between the invariants depending on two contributions: one is global, and it is related to the properties of the general fibre (the number $\lambda$), the other is local, and it depends on the geometry of the surface near the special fibres (the Horikawa number of the atoms).\\

The existence of a slope equality for a prescribed general condition is still an open question. It is not known, for example, if there is a slope equality for fibrations whose general fibre belongs to an assigned subvariety of $\cal{M}_g$. Moreover, it is not known if a slope equality, if existing, is unique for that class of fibrations, i.e. if the number $\lambda$ and the Horikawa numbers are unique.\\

As pointed out by Reid in \cite{Reid2}, a very useful method to find relation of this kind is the analysis of the relative canonical algebra of a fibration $f\!: S \fr B$, i.e. the $\cal{O}_B$-algebra of $\cal{R}(f)=\bigoplus_{n \geq 0} \cal{R}_n$ where the sheaves $\cal{R}_n$ are defined as \mbox{$\cal{R}_n=f_* \omega_{f}^{\otimes n}$}.\\

One of the core results of this theory is due to Horikawa (\cite{Hori}, 1977) which, in the case of genus 2 fibrations, established the following:
\begin{teo}[Horikawa]
Let $f\!: S \fr B$ a surface fibration with fibres of genus 2. Then
\begin{equation}
K_f^2
= 2\chi_f  +\sum_{P \in B}
\cal{H}(f,F_P)
\end{equation}
where the Horikawa number of a genus 2 fibre germ over $P$ 
\[ \cal{H}(f,F_P)=\operatorname{length} \operatorname{coker}((\text{\emph{Sym}}^2\cal{R}_1)_P \fr (\cal{R}_2)_P)\]
can be interpreted (roughly speaking) as the virtual number of  $\,2$-disconnected fibres of type $E_1 + E_2$ (with
$E_1$, $E_2$ elliptic curves meeting transversally in one point).\end{teo}

In her PhD thesis Mendes Lopes completed the local analysis of the canonical algebra of curves of genus 2 and 3 (see \cite{ML}). Based on this, a work of Reid (\cite{Reid2}, 1990) led to a proof of the following statement:

\begin{teo}[Reid]\label{Reid}

Let $f\!: S \fr B$ a surface fibration with fibres of genus 3 and suppose that the general fibre is nonhyperelliptic. Then
\begin{equation}
K_f^2
= 3\chi_f  +\sum_{P \in B}
\cal{H}(f,F_P)
\end{equation}
where the Horikawa number $\cal{H}(f,F_P)$ is defined as \[\cal{H}(f,F_P)=\operatorname{length} (\operatorname{coker}((\text{\emph{Sym}}^2\cal{R}_1)_P \fr (\cal{R}_2)_P)\emph{)}.\]
\end{teo}

A further, more general development is due to Konno (see \cite{K}, 1999 and \cite{AK}, 2002). He proves a slope equality for surface fibrations with fibre of odd genus $g$, with the general condition that the generic fibre has maximal gonality. He shows that there is the following relation
\begin{equation}K_f^2= \frac{6 (g-1)}{g+1} \chi_f + \sum_{P \in B} \cal{H}(f, F_P)\end{equation}
where the Horikawa indexes are calculated by means of the Koszul complex
\[\begin{matrix}
0 \fr \Lambda^{(g-1)/2+1} \cal{R}_1 \fr \ldots \fr \Lambda^{i} \cal{R}_1 \otimes \cal{R}_{(g-1)/2+1-i} \fr \Lambda^{i-1} \cal{R}_1 \otimes \cal{R}_{(g-1)/2+2-i} \fr \\ 
\ldots \fr \cal{R}_1 \otimes \cal{R}_{(g-1)/2} \fr \cal{R}_{(g-1)/2+1}\fr 0.
\end{matrix}\]
In particular he shows that a smooth fibre with maximal gonality is not an atom. In the case of genus 3 the work of Konno coincide precisely with the results of Reid as stated in Theorem \ref{Reid}. \\

The importance of the gonality in these works is based on the famous Green's conjecture (\cite{green}, 1984) which relates the cohomology of the Koszul complex of a smooth curve to its gonality. In particular it says that the Koszul complex is exact over a curve of maximal gonality.\\

Not many other slope equalities are known at the moment. One is due again to Ashikaga and Konno (see \cite{AK}, \S 2.2), which study hyperelliptic fibrations, i.e. the general condition they ask is that the generic fibre is smooth and hyperelliptic.\\

In this paper we are interested in a explicit slope equality for fibration of genus 5 with general fibre of maximal gonality. In fact, even thought the result of Konno in \cite{AK} and \cite{K} is very general and deep, it does not give an explicit and geometric explanation of the meaning of the Horikawa numbers, nor it describes explicitly the nature of the atoms. 

We show that, under the additional condition that every fibre is canonical, the same slope equality found by Konno for fibration of genus 5 can be computed in a straightforward and geometric way. 

The key of our result is the proof of Green's conjecture in the case of any singular canonical curve. In the genus 5 case Green's conjecture is equivalent to Noether and Petri theorems, thus our result becomes:\\

\noindent \textbf{Proposition \ref{lemma}} \emph {Let $C$ be a Gorenstein, genus 5 curve. Let us suppose that $C$ is $3$-connected and nonhyperelliptic (in the sense of Definition \ref{def hyperl}).} 

\emph{Then the dualizing sheaf $\omega_C$ is very ample. Moreover Noether's Theorem holds, i.e.
\begin{equation*}\Sym^n H^0(C, \omega_C) \fr H^0(C, \omega_C^{\otimes n})\end{equation*}
is surjective for any $n \geq 1$.}

\emph{Also Petri's Theorem holds. In the case of a genus 5 curve it means that the ideal of its canonical embedding in $\mathbb{P}^4$ is defined by three quadrics if the curve is nontrigonal, and by three quadrics and two cubics if it is trigonal.}\\

This guarantees that the contribution of a fibre depends only on its gonality and not on the kind of singularity that the fibre has (where the gonality for a genus 5 singular curve is defined in Definition \ref{def trig} and \ref{def hyperl}). Our main result is the following:\\

\noindent \textbf{Theorem \ref{enunciato}} \emph {Let $S$ be a projective surface with at most Du Val singularities and $B$ a projective smooth curve. Let $f\!: S \fr B$ be a
surface fibration of genus $5$ and let us suppose that every fibre is $3$-connected and not honestly hyperelliptic. Let us also suppose that the general fibre is nontrigonal.} 

\emph{Then there is a skyscraper sheaf $\cal{F}$ over $B$ such that if we define \[ \cal{H}(f, F_b)=1/2 \operatorname{length}(\cal{F}_b)\] we have that}
\begin{equation*}
K_f^2=4 \chi_f + \sum_{b \in B} \cal{H}(f, F_b)
\end{equation*}

\emph{Moreover $\operatorname{length}(\cal{F}_b)$ is even and strictly positive for any $b \in B$ supporting a trigonal fibre, and zero elsewhere.}\\

We compared this slope equality to the one found by Konno in \cite{K} and \cite{AK}. Using the notation introduced in Section \ref{section Konno}, in the genus 5 case it can be stated in the following way:
\begin{flushleft} \textbf{Theorem \ref{Konno} (Konno)} \emph{Let $f\!: S \fr B$ be a surface fibration with fibres of genus 5. Suppose moreover that the gonality of the general fibres is maximal. Then  we have that}\end{flushleft}
\[K_f^2 = 4 \chi_f + \sum_{P \in B} \frac12 (\operatorname{length} (\cal{K}_{1,2})_P-\operatorname{length} (\cal{K}_{0,3})_P).\]
\emph{Moreover $\operatorname{length} (\cal{K}_{1,2})_P-\operatorname{length} (\cal{K}_{0,3})_P$ is non-negative for any $P \in B$ and it vanishes on any smooth fibre of maximal gonality.}\\

A comparison of the two slope equality lead to the fact that they behave properly, i.e. they coincide. This means that the local analysis we did in order to prove Theorem \ref{enunciato} enlightens the geometric nature of the Horikawa index in Konno's relation.\\

\noindent \textbf{Proposition \ref{confronto}} \emph { Let $f\!: S \fr B$ a fibration satisfying the assumptions of Theorem \ref{enunciato}. Then the Horikawa numbers of Theorem \ref{enunciato} and of Theorem \ref{Konno} coincide. More precisely, the sheaf $\cal{K}_{0,3}$ is everywhere 0 and the sheaf  $\cal{K}_{1,2}$ is supported on the points with a trigonal fibre and is isomorphic to the sheaf $\cal{F}$, defined in Theorem \ref{enunciato} (that ``counts'' the trigonal fibres).}\\

Since the proof of Theorem \ref{enunciato} clarifies the nature and the structure of the atoms, we are able to construct many examples of regular surfaces satisfying the assumptions of Theorem \ref{enunciato} for every possible value of $p_g=h^0(S, \omega_S)$. The idea is based on the fact that any surface of this kind is a relative canonical model, thus $S= \textbf{Proj}_B(\cal{R}_f)$ and it has an embedding over $B$ in the $\mathbb{P}^4$-bundle $\textbf{Proj}_B(\bigoplus_{n \geq 0} \sym^n f_* \omega_f)$, which restricts to the canonical embedding of every curve. In Section 3 we showed that the ideals of these embeddings have a very specific form, i.e. they are generated by Pfaffian equations. In Section 6 we exhibit explicit $\mathbb{P}^4$-bundles over $\p$ and explicit equations in order to find surfaces with assigned invariants and atoms.  \\

\noindent \textbf{Acknowledgments:} I would like to thank my supervisor Miles Reid for introducing me to this beautiful topic and for all the help he generously gave me. I am deeply in debt to Lidia Stoppino for many suggestions, constant support and encouragement. I am grateful to the geometry groups at Pavia and Warwick University for their help.

\section{The relative canonical algebra of a fibration}
In this section we recall some standard facts and notations related to the relative canonical algebra of a fibred surface (see \cite{AK} and \cite{Reid2}).\\

A \textit{surface fibration} is a morphism between a connected projective surface with at most Du Val singularities and a smooth projective curve over an algebraically closed field $k$ of characteristic 0. A fibration which has no $-1$-curves in any of its fibres is called \textit{relatively
minimal}. We will always assume that this is the case.

This morphism is always flat (see \cite{Hart}, Proposition III.9.7).
For $b \in B$ we define $F_b=f^{*}(b)$ for a fibre of $f$, while $F$ is a general fibre. Since char$(k)=0$ Ramanujan's Lemma (see \cite{BPV}, Proposition III 11.1) implies that $h^0(F_b,\cal{O}_{F_b})=1$ for any $b \in B$. Thus $h^0(F_b, \omega_{F_b})$ is constant and we call this common value $g$, the genus of a smooth fibre.

Let $\omega_f$ be the \emph{relative canonical bundle} given by \[\omega_f=\omega_S \otimes f^*\omega_B^{-1}\] and $K_f$ the associated divisor. It is clear by adjuction that ${\omega_f}_{|F_b}\cong \omega_{F_b}$ for any $b \in B$. Let \[\chi_f=\chi(\cal{O}_S)- \chi(\cal{O}_B)\,\chi(\cal{O}_F)\] be the \emph{relative Euler-Poincar\'e characteristic.}

\begin{defi} We define the {\em relative canonical model} of a surface fibration \mbox{$f\!: S \fr B$} as the minimal surface fibration $f'\!: X \fr B$ such that $X$ is birationally equivalent to $S$ over $B$ and $\omega_{f'}$ is relatively ample.
In particular $X$ is constructed by contracting all the $-1$ and $-2$-curves of $S$ contained in the fibres of $f$.
\end{defi}

The fundamental invariants associated to a surface fibration $f\!: S \fr B$ are:

\begin{enumerate}
\item the self-intersection $K_f^2$ of the relative canonical divisor,
\item the relative Euler characteristic $\chi_f= \chi (\cal{O}_S)- \chi(\cal{O}_B)\,\chi(\cal{O}_F)$,
\item the relative topological Euler characteristic \[e_f=e(S) - e(F)\,e(B)=\sum_{b \in B}( e(F_b)-2+2g)\] where $e(.)$ is the topological Euler number of a space.
\end{enumerate}

The three invariants are related by Noether's Formula (see \cite{Bea}, Proposition~I.14):
\begin{equation}
K_f^2+e_f=12 \chi_f,
\end{equation}
thus only two of them are independent.
The invariants of a fibration are equal to the one of its relative canonical model $f'$, since the pullback of $\omega_{f'}$ is precisely~$\omega_f$.

\begin{defi}The \textit{relative canonical algebra} of the fibration $f$ is defined as 
\begin{equation}
\cal{R}(f)=\bigoplus_{n \geq 0} \cal{R}_n
\end{equation}
where the sheaves $\cal{R}_n$ are defined as
\begin{equation}
\cal{R}_n=f_* \omega_{f}^{\otimes n}.
\end{equation}
with multiplication induced by the tensor product $\omega_{f}^{\otimes n} \otimes \omega_{f}^{\otimes m}
\fr \omega_{f}^{\otimes n+m}$. It is a finitely generated $\cal{O}_B$-algebra, generated in degree
$\leq 4$ (the 1-2-3 Theorem, see \cite{K}).\end{defi}

The fibres of $\cal{R}_n$
\[\cal{R}_n \otimes k(b)=H^0(F_b,{(\omega_{f}^{\otimes n}})_{|F_b} )=H^0(F_b, \omega_{F_b}^{\otimes n})\]
have constant dimension, so by base change (see \cite{Mum}, Corollary II.2), $\cal{R}_n$ is a locally free sheaf. The rank of $\cal{R}_n$ is given by

\begin{equation}\label{rank}\text{rank}(\cal{R}_n)=\left\{ \begin{array}{ll}
1 & n=0\\
g & n=1\\
(2n-1)(g-1) & n >1.\\
\end{array}\right.\end{equation}

We can easily calculate $\chi(\cal{R}_n)$. The Leray spectral sequence gives \[H^n(S, \cal{F})=\bigoplus_{p+q=n}H^p(B, R^qf_*\cal{F})\] for any coherent sheaf $\cal{F}$ on
$S$. Again by base change we find that \[R^1f_*\omega_{f}^{\otimes n}=0 \text{ for all }n>1\] and
\[R^2f_*\omega_{f}^{\otimes n}=0 \text{ for all }n \geq 1.\] Moreover Grothendieck duality implies that \[R^1f_* \omega_{f}=\fs_B.\]
Thanks to these relations we find

 \begin{equation} \left\{ \begin{array}{ll}

\chi(\cal{R}_1)=\chi(\omega_{f})-\chi(\fs_B)&\\
\chi(\cal{R}_n)=\chi(\omega_{f}^{\otimes n})& \text{for all  } \, n >1.

\end{array}\right.\end{equation}

Using the Riemann-Roch formula for surfaces and curves we finally compute
\begin{eqnarray} \label{chi}
\chi(\cal{R}_n) &=&\left\{ \begin{array}{ll}
b-1 & n=0,\\
\chi_f+g(1-b)& n=1,\\
\chi_f + {n \choose 2} K_f^2-(2n-1)\,\chi(\cal{O}_B)\,\chi(\cal{O}_F)& n >1;\\
\end{array}\right.\\[6pt]\label{deg}
 \deg(\cal{R}_n) &=&\left\{ \begin{array}{ll}
0 & n=0,\\
\chi_f& n=1,\\
\chi_f + {n \choose 2} K_f^2& n >1.\\
\end{array}\right. \end{eqnarray}

Since we are interested in the analysis of the relative canonical algebra, it is important to understand its local structure, i.e. the canonical algebra of a curve. There is a deep relation between the properties of the canonical sheaf of a Gorenstein curve and some connection properties.

\begin{defi}
A Gorenstein curve $C$ is said to be \textit{numerically $k$-connected} if
\begin{equation}
\deg (\omega_C)_{|B}- \deg \omega_B  \geq k
\end{equation}for every generically Gorenstein strict subcurve $B \subset C$.

If $C$ is an effective divisor on a smooth surface, it is numerically $k$-connected if and only if, whenever we write $C = A+B$ as a sum
of effective divisors $A, \, B $, we have that $A.B \geq k$. In this case the sheaf $(\omega_C)_{|B} \cong \omega_B \otimes \cal{O}_B(A)$ (see \cite{BPV} II.6) and its degree is precisely $\deg \omega_B+ A.B$.
\end{defi}

\begin{defi}\label{def trig} A smooth curve is \textit{trigonal} if it has a $g_3^1$ without base points.
A Gorenstein curve is \textit{trigonal} if it is a limit of trigonal smooth curves but not of hyperelliptic ones, in the sense that the curve is a fibre of a flat family of curves such that the generic fibre is smooth and trigonal, but there is no such family with hyperelliptic general fibre. 
\end{defi}
\begin{rmk}
The locus of trigonal curves has codimension 1 in the moduli space of genus 5 stable curves $\overline{\cal{M}}_5$, while the general curve has gonality 4.
\end{rmk}

\begin{defi}[\cite{Ca}, Definition 3.18]\label{def hyperl}
We say that a Gorenstein curve $C$ is {\em honestly hyperelliptic} if there exists a finite morphism  $\psi\!:C \fr \p$ of degree 2 (that is, $\psi$  is finite and  $\psi_* \cal{O}_C$ is locally free of rank 2 on $\p$).
The linear system  $\psi^*|\cal{O}_{\p}(1)|$ defining $\psi$  is called an honest $g_2^1$.
\end{defi}

The key theorem we are going to need is the following:

\begin{teo}[\cite{CFHR}] Let $C$ be a numerically 3-connected Gorenstein curve. Then either $\omega_C$ is very ample or $C$ is honestly hyperelliptic. In particular if $p_a(C) \geq 2$ then $\omega_C$ is ample, and if $p_a(C)=1$ then $C$ is honestly hyperelliptic.
\end{teo}
 
\section{Fibrations of genus 5}
In this section we present the proof of Theorem \ref{enunciato}. Since the proof is based on a global version of the analysis of the canonical map of the fibres, it is essential to understand the nature of such maps. The following proposition shows that, for our purposes, a canonical singular curve is by no means different from a nonhyperelliptic smooth one. 

\begin{prop}\label{lemma}
Let $C$ be a Gorenstein, genus 5 curve. Let us suppose that $C$ is $3$-connected and nonhyperelliptic (in the sense of Definition \ref{def hyperl}).\\

Then the dualizing sheaf $\omega_C$ is very ample. Moreover Noether's Theorem holds, i.e.
\begin{equation}\label{noether genus 5}\Sym^n H^0(C, \omega_C) \fr H^0(C, \omega_C^{\otimes n})\end{equation}
is surjective for any $n \geq 1$.

Also Petri's Theorem holds. In the case of a genus 5 curve it means that the ideal of its canonical embedding in $\mathbb{P}^4$ is defined by three quadrics if the curve is nontrigonal, and by three quadrics and two cubics if it is trigonal.
\end{prop}

\begin{proof}
Thanks to Theorem 3.6 in \cite{CFHR}, we know that the dualizing sheaf $\omega_C$ is very ample precisely when the curve is 3-connected and not honestly hyperelliptic.\\

It is already known that the map (\ref{noether genus 5}) is surjective for $n \geq 3$ (see \cite{K}, Proposition 1.3.3) whenever it is surjective for $n=2$. \\

It remains to be proved that there can not be more that 3 linearly independent quadrics in $\mathbb{P}^4$ vanishing on the curve. There are at least 3 quadrics, since a smooth curve is cut out by exactly 3 quadrics (Petri's Theorem, see \cite{ACGH}, Chapter III). Thus we consider at first the intersection of any 3 of these quadrics. We distinguish two cases.

\begin{enumerate}
\item The intersection is a curve: then it is a degree 8 curve since it is a complete intersection of 3 quadrics, and it must coincide with $C$. Then the ideal is generated by exactly 3 quadrics.
\item The intersection is a surface $S$: $S$ must be non degenerate since it contains $C$, so it must have degree at least 3, since no degree 1 or 2 surface can span the whole $\mathbb{P}^4$. But it is a strict subset of the intersection of 2 quadrics, which is a surface of degree 4, so $S$ has exactly degree 3. If there is a fourth quadric, the intersection between the surface and the quadric can be either a smaller surface, impossible again since the degree can not drop, or a curve of degree 6, which must contain $C$ of degree 8, impossible again, so there is no fourth quadric.
\end{enumerate}

We have proved that there is exactly a 3 dimensional vector space of quadrics in the ideal of the curve, but
\[3=\dim \Sym^2H^0(C, \omega_C)-h^0(C, \omega_C^2)\]
so Noether's Theorem holds.
We have also shown that the intersection of the 3 quadrics in the ideal of the curve can be either the curve itself, or a degree 3 surface $S$. The surface S is not necessarily irreducible. Again we distinguish two cases.
 \begin{enumerate}
\item The surface $S$ is irreducible: it can only be a rational normal scroll of the form $\mathbb{F}(1,2)$ or $\mathbb{F}(0,3)$ embedded with $\cal{O}_{\mathbb{F}}(0,1)$ (see \cite{Reid1}, Section 2.12). Thanks to a straightforward application of the adjunction formula one can see that the curve is defined, in the scroll, by a unique equation of relative degree 3.
\item The surface $S$ is reducible: the union of a degree 1 surface ($\mathbb{P}^2$) and a quadric surface, which could be reducible itself. Since the surface is not degenerate, and the curve is 3-connected (since $\omega_C$ is very ample), the intersection of the 2 components must be a line. With straightforward calculations one can prove that the equations of the surface are in the form
\begin{equation}
\operatorname{rank}\left( \begin{matrix} x_0 & l(x) & m(x)\\ 0 & x_1 & x_2\end{matrix}\right)\leq 1
\end{equation}
\end{enumerate}
In both cases one can check that there are exactly 2 linear syzygies between the quadrics in the ideal of the curve. But there is a 15 dimensional vector space of cubics in the ideal (thanks to Riemann-Roch Theorem), a 3-dimensional vector space of quadrics, 5 linear forms $x_i$ and 2 syzygies, so there are exactly 2 more cubics that we need to add to the quadrics in order to cut the curve.\\

We will see now that a curve is an intersection of quadrics if and only if it is nontrigonal. If a curve is a flat limit of smooth trigonal curves, it is clear that its ideal will need at least the same number of generators as the case of the general curve, so the 2 additional cubics are needed. Conversely, if a curve is the intersection of a degree 3 surface of the kind we described with 2 cubic equations, it is not difficult, using Bertini's Theorem, to modify its equations in a smooth way in order to build a family of smooth curves, which are intersection of a smooth scroll $\mathbb{F}(1,2)$ and two cubic equations, thereby trigonal curves.
\end{proof}

We are now in a position to prove our main result.
\begin{teo}\label{enunciato}

Let $S$ be a projective surface with at most Du Val singularities and $B$ a projective smooth curve. Let $f\!: S \fr B$ be a
surface fibration of genus $5$ and let us suppose that every fibre is $3$-connected and not honestly hyperelliptic. Let us also suppose that the general fibre is nontrigonal.

Then there is a skyscraper sheaf $\cal{F}$ over $B$ such that if we define \[ \cal{H}(f, F_b)=1/2 \operatorname{length}(\cal{F}_b)\] we have that

\begin{equation}\label{teo}
K_f^2=4 \chi_f + \sum_{b \in B} \cal{H}(f, F_b)
\end{equation}

Moreover $\operatorname{length}(\cal{F}_b)$ is even and strictly positive for any $b \in B$ supporting a trigonal fibre, and zero elsewhere.

\end{teo}

\begin{proof} Consider the map \[\phi_n\!: \text{Sym}^n(\cal{R}_1) \fr \cal{R}_n\] Its fibre over a point $t \in B$ is just

\[\sym^n H^0(F_t, \omega_{F_t}) \fr H^0({F_t}, \omega_{F_t}^{\otimes n}).\]
By Noether's Theorem this map of stalks is surjective, thus $\phi_n$ is surjective.

We are interested in the locally free sheaves $\cal{K}_n= \text{ker}( \phi_n)$. For each $t \in
B$, the fibre $\cal{K}_n \otimes k(t)$ is the vector space of polynomials of degree $n$ in $\mathbb{P}^4$
vanishing on the curve $F_t$:
\begin{equation}
\cal{K}_n \otimes k(t)= H^0(\mathbb{P}^4, I_{F_t}(n)).
\end{equation}
In particular one can see that rank$(\cal{K}_2)=3$ and rank$(\cal{K}_3)=15$.\\

Let us consider the commutative diagram:

\begin{equation}\label{diagramma}\xymatrix{
0 \ar[r]&  \cal{K}_2\otimes \cal{R}_1 \ar[r] \ar[d]^{\mu}& \text{Sym}^2 \cal{R}_1 \otimes \cal{R}_1 \ar[r] \ar[d]& \cal{R}_2 \otimes \cal{R}_1 \ar[r]\ar[d]  & 0\\
0 \ar[r] & \cal{K}_3 \ar[r] & \text{Sym}^3 \cal{R}_1 \ar[r]& \cal{R}_3 \ar[r]& 0 }
\end{equation}

The fibre of the map $\mu\!: \cal{K}_2 \otimes \cal{R}_1 \fr \cal{K}_3  $ over $t \in B$ is 
\[\mu_t\!: H^0(\mathbb{P}^4, I_{F_t} (2))\otimes H^0(\mathbb{P}^4,\cal{O}_{\mathbb{P}^4}(1))  \fr H^0(\mathbb{P}^4, I_{F_t} (3)).  \]
We can apply Lemma \ref{lemma} and see that for any point $t$ supporting a nontrigonal fibre, $\mu_t$ is surjective, since the ideal $I_{F_t}$ is generated by quadrics. In particular $\mu_t$ is an isomorphism since the two vector spaces have both dimension 15 over the base field $k$. Conversely, over a point supporting a trigonal fibre, the cokernel of $\mu_t$ is a 2 dimensional vector space.

The kernel of the map of sheaves $\mu$ is trivial, while the cokernel is a skyscraper sheaf $\cal{F}$ supported on the points with trigonal fibres:

\begin{equation}\label{quad e cubiche}
0 \fr \cal{K}_2 \otimes \cal{R}_1 \fr \cal{K}_3 \fr \cal{F} \fr 0.
\end{equation}

Taking into account formulas (\ref{chi}) and (\ref{deg}) it is possible to calculate the Euler characteristic of $\text{Sym}^2(\cal{R}_1)\otimes \cal{R}_1$ and $\text{Sym}^3(\cal{R}_1)$ using the splitting principle:
\begin{eqnarray}
\deg (\text{Sym}^2(\cal{R}_1))&=& 6 \chi_f, \\
\deg (\text{Sym}^3(\cal{R}_1))&=& 21 \chi_f. \nonumber
\end{eqnarray}

Using standard results on Chern classes of a tensor product, one can calculate that
\begin{eqnarray}
\deg (\text{Sym}^2(\cal{R}_1)\otimes \cal{R}_1)&=& 45 \chi_f \\
\deg (\cal{R}_2\otimes \cal{R}_1)&=& 17 \chi_f+5 K_f^2 \nonumber
\end{eqnarray}
and if we take into account diagram (\ref{diagramma}) we can conclude that
\begin{eqnarray}\label{chi strani}
\chi(\cal{K}_3)&=&20 \chi_f-3 K_f^2 +15 \chi(\cal{O}_B)\\
\chi(\cal{K}_2\otimes \cal{R}_1)&=&28 \chi_f-5 K_f^2 +15 \chi(\cal{O}_B).\nonumber
\end{eqnarray}

Now consider the exact sequence (\ref{quad e cubiche}). By additivity of the Euler characteristic and substituting the equality (\ref{chi strani}) we find
\begin{eqnarray*} 0 &= &\chi (\cal{R}_2 \otimes \cal{R}_1)- \chi (\cal{R}_3)+ \chi (\cal{F})\\
&=& (28 \chi_f-5 K_f^2 +15 \chi(\cal{O}_B))-(20 \chi_f-3 K_f^2 +15 \chi(\cal{O}_B))+\chi (\cal{F})\\
&=& 8 \chi_f-2 K_f^2 + \chi (\cal{F}).
\end{eqnarray*}
This is precisely equation (\ref{teo}) of Theorem \ref{enunciato}.
\end{proof}
\newpage
As we have seen, the sheaf $\cal{F}$ represents the extra cubics needed to define the ideal of a trigonal curve in its canonical embedding. In the next proposition we will study further the structure of this sheaf, and in particular we will prove that $\operatorname{length} (\cal{F}_t)$ is even for any $t \in B$. 

\begin{prop}\label{local struct} Let $f\!: S \fr B$ a surface fibration satisfying the assumptions of Theorem \ref{enunciato}. Let us fix a point $b \in B$ supporting a trigonal fibre and $U$ a small neighborhood of $b$, with no other trigonal fibres.

Then the equations of $S_U=f^{-1}(U) \subset U \times \mathbb{P}^4$ are the Pfaffians of the following skew matrix with the upper triangular entries as follows:
\begin{equation}M=\left( \begin{matrix}
  & t^n & l_1(x) & l_2(x) & l_3(x)\\
 & & m_1(x) & m_2(x) & m_3(x)\\
 & & & q_3^t(t,x) & -q_2^t(t,x)\\
 & & & & q_1^t(t,x)
\end{matrix} \right)
\end{equation}
where $t$ is a local parameter in $U$ centered in $b$, $x=(x_0, \ldots, x_4)$ are the variables in $\mathbb{P}^4$, $l_i(x)$ and $m_i(x)$ are linear and $q_i^t(t,x)$ are polynomials of degree 2 in the variables $x$.
\end{prop}
\begin{rmk}\label{equaz} For those unfamiliar with the Pfaffian notation, the Pfaffian equations of $M$ are:
\begin{equation}\label{equaz locali} \operatorname{Pf}(M)=\left \{\begin{array}{lll}
c_1(t,x)&=&\sum_i m_i(x)q_i^t(t,x) \\
c_2(t,x)&=&\sum_i l_i(x)q_i^t(t,x)\\
p_1(t,x)&=&(m_2(x)l_3(x)-m_3(x)l_2(x))+t^n q_1^t(t,x)\\
p_2(t,x)&=&(m_1(x)l_3(x)-m_3(x)l_1(x))+t^n q_2^t(t,x)\\
p_3(t,x)&=&(m_1(x)l_2(x)-m_2(x)l_1(x))+t^n q_3^t(t,x)
\end{array}\right.\end{equation}
\end{rmk}

\begin{rmk} We prove with a local analysis that the equation are in the Pfaffian form of a $5 \times 5$ matrix. Anyway a more general result by Buchsbaum and Eisenbud in \cite{BE} garantees that any codimension 3 ideal must be in the Pfaffian form.
\end{rmk}

\begin{proof}  Over $U$, Petri's Theorem tells us that $S_U$ is the intersection of three quadrics $q_1(t,x),\,q_2(t,x),\,q_3(t,x)$ plus two cubics $c_1(t,x)$ and $c_2(t,x)$. The cubics are in the ideal of every fibre except for the one over $t=0$, thus
\begin{eqnarray*}
c_1(t,x)&=& \sum m_{i}(t,x)q_i(t,x)\\
c_2(t,x)&=& \sum l_i(t,x)q_i(t,x)
\end{eqnarray*}
where $m_i(t,x)$ and $l_i(t,x)$ are linear in $x$ and well defined away from $t=0$, or, equivalently:
\begin{eqnarray}\label{ni}
t^{n_1}c_1(t,x)&=& \sum m_{i}(t,x)q_i(t,x)\\
t^{n_2}c_2(t,x)&=& \sum l_i(t,x)q_i(t,x)\nonumber
\end{eqnarray}
with $m_i(t,x)$ and $l_i(t,x)$ regular on the whole $U\times \mathbb{P}^4$.

By changing the $c_i(t,x)$ without modifying the surface, we can suppose that $m_i(t,x)$ and $l_i(t,x)$ depend only on $x$.\\

Write
\begin{equation}
q_i(t,x)=q_i^0(x)+t^{d_i}q_i^t(t,x)
\end{equation}
with $d_i$ minimal. Equations (\ref{ni}) become:
\begin{eqnarray}\label{n1 n2}
t^{n_1}c_1(t,x)&=&\sum m_{i}(x)q_i^0(x) +\sum m_i(x)t^{d_i}q_i^t(t,x)\\
t^{n_2}c_2(t,x)&=&\sum l_{i}(x)q_i^0(x) +\sum l_i(x)t^{d_i}q_i^t(t,x). \nonumber
\end{eqnarray}

Thus on $\mathbb{P}^4$ we have that \begin{eqnarray}\label{syzy}\sum a_{i}(x)q_i^0(x)=0,\\
\sum b_{i}(x)q_i^0(x)=0.\nonumber
\end{eqnarray}

We can check that, up to scalar multiplication, we have that
\begin{equation}\label{defq^0} \left \{\begin{array}{lll}
q_1^0(x)&=&m_2(x)l_3(x)-m_3(x)l_2(x)\\
q_2^0(x)&=&-m_1(x)l_3(x)+m_3(x)l_1(x)\\
q_3^0(x)&=&m_1(x)l_2(x)-m_2(x)l_1(x)
\end{array}\right.\end{equation}

Consider again equation (\ref{n1 n2}): the second term $\sum m_i(x)t^{d_i}q_i^t(t,x)$ is not divided by $t^{n_1}+1$, and similarly for $c_2$.

Hence there is one $d_i$ (say $d_1$) equal to $n_1$ and the same for $c_2$. If $n_1=n_2$ we have concluded the proof: the 5 generators of the ideal $q_0^1,\, q_0^2,\,q_0^3,\,c_1, c_2$ are in the required form.\\

We want to show that the converse is impossible. Let us suppose that  $n_1$ and $n_2$ are different, say $n_1<n_2$ and $d_2=n_2$. Taking into account equation (\ref{syzy}), equation (\ref{n1 n2}) becomes
\begin{eqnarray}
c_1(t,x)&=& m_1(x)q_1^t(t,x)+ t^{n_2-n_1 }m_2(x)q_2^t(t,x)+t^{d_3-n_1}m_3(x)q_3^t(t,x)\nonumber\\
c_2(t,x)&=&l_2(x) q_2^t(t,x)+t^{-n_2}(t^{n_1}l_1(x)q_1^t(t,x)+t^{d_3}l_3(x)q_3^t(t,x))
\end{eqnarray}

This implies that $d_3=n_1$ and that
\begin{equation}\label{zero in pp4}l_1(x) q_1^t(0,x)+l_3(x)q_3^t(0,x)=0 \text{  on } \mathbb{P}^4.\end{equation}

We have also that
\begin{eqnarray}\label{ci}
c_1(0,x)&=& m_1(x)q_1^t(0,x)+m_3(x)q_3^t(0,x)\\
c_2(0,x)&=&l_2(x) q_2^t(0,x)\nonumber
\end{eqnarray}
A consequence of the last two equations is that on $\mathbb{P}^4$
\begin{eqnarray}
l_1(x)c_1(0,x)&=& q_2^0(x)q_3^t(0,x)\\
l_3(x)c_1(0,x)&=& q_2^0(x)q_1^t(0,x).\nonumber
\end{eqnarray}
Thus one of the following holds:
\begin{enumerate}
\item $l_1=l_3=0$ on $\mathbb{P}^4$. This is impossible since this would imply $q_2^0=0$ on $\mathbb{P}^4$ because of equation (\ref{defq^0}), but $q_2^0$ is one of the generators of the ideal;
\item $q_2^0(x)$ divides $c_1(0,t)$. This is impossible because of the definition of $c_1$;
\item $l_1(x)$ and $l_3(x)$ divide $q_2^0$ and they are not multiple one of the other. Thus $q_2^0(x)= \beta l_1(x)l_3(x)$, but thanks to the last equation $c_1(0,x)$ is divided by $l_1(x)$ and $l_3(x)$, thus by $q_2^0(x)$. We already said this is impossible.
\item $l_3(x) = \beta l_1(x)$ and $l_1(x)$ divides $q_2^0(x)$ (or $l_1(x)=\beta l_3(x)$ etc.). In this case one can check that there is an irreducible component of the fibre over 0, defined by $l_1(x)=l_2(x)=q_3^t(x)=0$, which is a $-2$-curve, but we already showed that this can not be because every fibre is $3$-connected.
\end{enumerate}
Since every case is impossible, our assumption $n_1 < n_3$ is impossible.  \end{proof}

\begin{cor}\label{cor n} With the same notation of Proposition \ref{local struct} and Remark \ref{equaz} we have that
\begin{eqnarray*}
t^n c_1(t,x)&=&\sum m_i(x)p_i(t,x) \\
t^n c_2(t,x)&=&\sum l_i(x)p_i(t,x).
\end{eqnarray*}
In particular, the stalk of $\cal{F}$ in $b$ has length $2n$, with $n \geq 1$.
\end{cor}
\begin{proof} It is immediate to check the relations using the explicit form of the equations involved in Remark \ref{equaz}.
\end{proof}

\begin{rmk} It is clear from Proposition \ref{local struct} and Corollary \ref{cor n} that the Horikawa number is positive for any trigonal fibre, but its exact value depends only on the embedding of the curve in the surface, and not on the intrinsic geometry of the curve itself. In Section 5 we will show several examples in which the same curve has any possible Horikawa number in different fibrations. It is moreover clear that the Horikawa number depends only on a neighborhood of the fibre and not on the entire surface. 

In a forthcoming paper in collaboration with L. Stoppino we are proving that, as one can expect, when the trigonal fibre is stable the Horikawa number has a modular meaning. More precisely, when the atom is stable (but the fibration does not need to be stable itself) we can compute its Horikawa number as an intersection number in $\overline{\cal{M}}_5$. This should answer the question made by Ashikaga in \cite{As}, Chapter I.4, which asks whether the slope equality found by Konno is related to a similar slope equality that exists for stable fibrations computed by means of the Harris Mumford formula (see \cite{HM}) for the trigonal divisor. 
\end{rmk}

\section{Konno's slope equality}\label{section Konno}
Consider the Koszul complex of a surface fibration with fibres of genus 5:

\[0 \fr \Lambda^3 \cal{R}_1 \fr \Lambda^2 \cal{R}_1 \otimes \cal{R}_1 \fr \cal{R}_1 \otimes \cal{R}_2 \fr \cal{R}_3 \fr 0\]

with maps \[\Lambda^i \cal{R}_1 \otimes \cal{R}_{3-i} \stackrel{d_{i,3-1}}{\longrightarrow} \Lambda^{i-1} \cal{R}_1 \otimes \cal{R}_{3-i+1}\]

Define $\cal{K}_{i,j}= \ker(d_{i,j})/\operatorname{Im}(d_{i+1, j-1})$. In \cite{AK} \S 2.3 and \cite{K_Cliff} Konno proves a slope equality for fibrations of maximal gonality that, in the case of a genus 5 fibration, can be stated in the following way:

\begin{teo}[Konno]\label{Konno}
Let $f\!: S \fr B$ be a surface fibration with fibres of genus 5. Suppose moreover that the gonality of the general fibres is maximal. Then  we have that
\[K_f^2 = 4 \chi_f + \sum_{P \in B} \frac12 (\operatorname{length} (\cal{K}_{1,2})_P-\operatorname{length} (\cal{K}_{0,3})_P).\]
Moreover $\operatorname{length} (\cal{K}_{1,2})_P-\operatorname{length} (\cal{K}_{0,3})_P$ is non-negative for any $P \in B$ and it vanishes on any smooth fibre of maximal gonality.

\end{teo}

The core of the proof is the study of the Koszul complex by means of Green's conjecture, which implies that the stalk of the sheaves $\cal{K}_{i,j}$ on a point supporting a smooth curve of maximal gonality are 0. It remains unclear the precise geometrical meaning of those sheaves outside the ``good" points. In particular, it is not clear which kind of points have nonzero stalk. Our geometrical analysis in Proposition \ref{lemma} will lead us to a deeper understanding of the meaning of these objects. In fact we proved that Green's conjecture hold for canonical singular curves as well, thus the only important information on a canonical fibre is whether this fibre is trigonal or not.

\begin{prop}\label{confronto} Let $f\!: S \fr B$ a fibration satisfying the assumptions of Theorem \ref{enunciato}. Then the Horikawa numbers of Theorem \ref{enunciato} and of Theorem \ref{Konno} coincide. More precisely, the sheaf $\cal{K}_{0,3}$ is everywhere 0 and the sheaf  $\cal{K}_{1,2}$ is supported on the points with a trigonal fibre and is isomorphic to the sheaf $\cal{F}$, defined in Theorem \ref{enunciato} (that ``counts'' the trigonal fibres).
\end{prop}
\begin{proof}By definition \[\cal{K}_{0,3}=\ker(d_{0,3})/\operatorname{Im}(d_{1, 2}) = \frac{\cal{R}_{3}}{\operatorname{Im}(\cal{R}_1 \otimes \cal{R}_{2} \stackrel{d_{1,2}}{\longrightarrow} \cal{R}_3)}. \] The homomorphism $d_{1,2}$ is surjective since, on the fibre over $b \in B$, it becomes the product map
\[H^0(F_b, \omega_{F_b}) \otimes H^0(F_b, \omega_{F_b}^{\otimes 2}) \fr H^0(F_b, \omega_{F_b}^{\otimes 3}) \]
which is surjective thanks to Noether's Theorem, as seen in Proposition \ref{lemma}.

Let us study $\cal{K}_{1,2}$ which by definition is
\[ \cal{K}_{1,2}= \ker(d_{1,2})/\operatorname{Im}(d_{2, 1}) = \frac{\ker(\cal{R}_1 \otimes \cal{R}_{2} \stackrel{d_{1,2}}{\longrightarrow} \cal{R}_3)}{\operatorname{Im}(\Lambda^2 \cal{R}_1 \otimes \cal{R}_{1} \stackrel{d_{2,1}}{\longrightarrow} \cal{R}_1 \otimes \cal{R}_{2})}. \]

The map $d_{2,1}$ factors through $(\Lambda^2 \cal{R}_1 \otimes \cal{R}_{1})/\Lambda^3 \cal{R}_1$ since the following sequence is exact for any vector space and any vector bundle $V$:
\[ 0 \fr \Lambda^3 V \fr \Lambda^2 V \otimes V \fr V \otimes \Sym^2 V \fr \Sym^3 V \fr 0. \]

We can consider the following exact diagram which is an enlarged version of diagram \ref{diagramma}:

\begin{equation}\xymatrix{ & & 0\ar[d] & 0\ar[d]\\
&0 \ar[d] \ar@{-->}[r]&  \frac{\Lambda^2 \cal{R}_1 \otimes \cal{R}_{1}}{\Lambda^3 \cal{R}_1} \ar[d]\ar[dr]^{d_{2,1}} \ar[r]^{d_{2,1}}& \ker(d_{1,2})\ar[d] \ar[r] & \cal{K}_{1,2} \ar[r] &0 \\
0 \ar[r] & \cal{R}_1 \otimes \cal{K}_2 \ar[d] \ar[r] & \cal{R}_1 \otimes \Sym^2 \cal{R}_1 \ar[d]\ar[r] & \cal{R}_1 \otimes \cal{R}_2\ar[d]^{d_{1,2}} \ar[r] & 0\\
0 \ar[r] &\cal{K}_3 \ar[d] \ar[r] & \Sym^3 \cal{R}_1 \ar[d]\ar[r]& \cal{R}_3\ar[r]\ar[d]&0\\
&\cal{F}\ar[d] & 0 & 0\\
& 0
}\end{equation}

With a simple diagram chase we can prove that the map \[d_{2,1}\!: (\Lambda^2 \cal{R}_1 \otimes \cal{R}_{1})/\Lambda^3 \cal{R}_1  \fr \ker(d_{1,2})\] is injective (and the dotted line in the diagram give rise to an exact row) and, more importantly, that its cokernel $\cal{K}_{1,2}$ is isomorphic to $\cal{F}$.

\end{proof}

\section{Examples}\label{Examples}

In this Section we exhibit examples of fibrations satisfying the assumptions of Theorem \ref{enunciato}.
In particular we are looking for fibrations $f\!: S \fr \p$ with $S$ smooth projective surface and a single trigonal fibre.

\begin{es}\label{ex} The easiest ambient space $\mathbb{P}$ we can work with is a rational
normal scroll, or even $\p_{t_0:t_1}\times \mathbb{P}^4_{x_0:\,\ldots\, :x_4}$.\end{es}

We already know the local form of the equations near a point supporting a trigonal fibre (see Proposition \ref{local struct}), we want to transform them in global equations. Consider the Pfaffian equations of the following skew
matrix:

\begin{equation}
M=\left( \begin{matrix}
  & t_1^n & x_0 & x_2 & x_3\\
 & & t_0^nx_1 & t_0^nx_3 & t_0^nx_4\\
 & & & q_1 & q_2\\
 & & & & q_3
\end{matrix} \right)
\end{equation}
 $q_1,\,q_2,\,q_3$ are generic quadratic homogeneous polynomials only in the $x$ variables. The Pfaffian equations of the matrix $M$ are:

\begin{equation} \label{pfaff}\operatorname{Pf}(M)=\left \{\begin{array}{lll}
c_1&=& t_0^n(x_1q_3-x_3q_2+x_4q_1)\\
c_2&=&x_0q_3-x_2q_2+x_3q_1\\
p_1&=&t_1^nq_3+t_0^n( - x_2x_4+x_3^2)\\
p_2&=&t_1^nq_2+t_0^n(-x_0x_4+x_1x_3)\\
p_3&=&t_1^nq_1+t_0^n(-x_0x_3+x_1x_2)

\end{array}\right.\end{equation}

Over each $(t_0, t_1) \neq (1,0)$ the two cubic polynomials $c_1$ and $c_2$ in (\ref{pfaff}) are linear combinations of
the three quadric polynomials:
\begin{eqnarray*}
t_1^n c_1&=& t_0^n(x_1 p_1 - x_2 p_2+x_4 p_3)\\
t_1^n c_2&=& x_0 p_1 - x_2 p_2+x_3 p_3.
\end{eqnarray*}

Over $(t_0, t_1) = (1,0)$ we impose a trigonal fibre. But we know
that a nonsingular trigonal curve of genus 5 in $\mathbb{P}^4$ is the intersection of three cubic polynomials
and a rational normal scroll $\mathbb{F}(1,2)$ (see \cite{Reid1}). The equations of $\mathbb{F}(1,2)\subseteq \mathbb{P}^4$ are given by
\begin{equation*}\text{rank} \left( \begin{matrix}
 x_0 & x_2 & x_3\\
 x_1 & x_3 & x_4\\
\end{matrix} \right)\leq 1.\end{equation*}
These equations coincide with $p_1, \, p_2,\, p_3$ when $t_1=0$. One can check that the equations $c_1$ and $c_2$ cut a trigonal curve inside the scroll $\mathbb{F}(1,2)$ (see \cite{Reid1}, Chapter 2).\\

The only issue now is to choose $q_1, q_2, q_3$ such that $S$ is a nonsingular surface. In order to do this we apply
Bertini's Theorem over the subset $(t_1 \neq 0) \subset \p \times \mathbb{P}^4$ and conclude that
the three quadric polynomials are general enough to have nonsingular intersection. Thus $S$ is smooth away from the trigonal fibre. Then we apply Bertini's Theorem again over $t_1=0$ to conclude that the trigonal fibre is smooth for a generic triple $q_1,\, q_2,\, q_3$. Thus for the generic triple $q_1,\, q_2,\, q_3$ $S$ is a nonsingular surface.\\

In this situation the projection $f\!: S \fr \p$ is a flat morphism, because it is a surjective
morphism between a smooth surface and a curve. Moreover, every fibre is connected since $h^0(F,
\cal{O}_F)$ is constant and we have chosen the fibre over $t_1=0$ to be a smooth connected curve. So every
fibre is a canonical genus 5 connected curve, thus we have no hyperelliptic fibres. They must be all
3 connected. In fact being 3 connected and nonhyperelliptic is equivalent to the existence of the canonical
embedding, and all the fibres are canonically embedded (see \cite{CFHR}).\\

Let us compute the invariants of the surface $S$. For this we need the free resolution of $S$ inside $\p \times \mathbb{P}^4$.

One can prove this is given by
\begin{equation}\label{risoluzione}
\begin{matrix}
0 \leftarrow \cal{O}_S \leftarrow \cal{L}_0 \stackrel{\text{Pf}
\,M}{\longleftarrow} \cal{L}_1
 \stackrel{M} {\leftarrow} \cal{L}_2 \stackrel{^t\text{Pf}
\,M}{\longleftarrow}\cal{L}_3\leftarrow 0\end{matrix}\end{equation}
where the $\cal{L}_i$ are direct sums of line bundles on $\p \times \mathbb{P}^3$ given by
\begin{eqnarray*}
\cal{L}_0&=&\cal{O}_{\p \times \mathbb{P}^4}\\
\cal{L}_1&=&\cal{O}_{\p \times \mathbb{P}^4}(-n,-3)\oplus \cal{O}_{\p \times
\mathbb{P}^4}(0,-3) \oplus 3\cal{O}_{\p \times \mathbb{P}^4}(-n,-2)\\
\cal{L}_2&=&\cal{O}_{\p \times \mathbb{P}^4}(-n,-3)\oplus \cal{O}_{\p \times
\mathbb{P}^4}(-2n,-3) \oplus 3\cal{O}_{\p \times \mathbb{P}^4}(-n,-4)\\
\cal{L}_3&=&\cal{O}_{\p \times \mathbb{P}^4}(-2n,-6)
\end{eqnarray*}
and $\operatorname{Pf} M$ is the row vector of Pfaffians of $M$.\\

Exactness can be checked on the fibres, thus one only needs to work out the free resolution of the canonical image of a nonhyperelliptic genus 5 curve. But this is straightforward both in the nontrigonal and in the trigonal case.\\

The dualizing sheaf $\omega_S$ can be computed by dualizing the free resolution of $\cal{O}_S$ (see \cite{Hart}, Proposition III.7.5). We conclude that
$\omega_S=\cal{O}_S(2n -2,1)$ and that there is an isomorphism  $H^0(\p \times \mathbb{P}^4,
\cal{O}(2n-2,1))\stackrel{\sim}{\fr}H^0(S,\omega_S)$.

In particular we see that $p_g=5(2n-1)$.\\

We can prove as well that the surface is regular. To prove this one considers again the resolution (\ref{risoluzione}) and check that most of the cohomology groups involved vanish.

The upshot of this is that $\chi_f=1+p_g+4=10n$.\\

It remains to compute $K_f^2$.

The relative canonical sheaf $\omega_{S|\p}=\omega_S\otimes f^*\cal{O}(2)=\cal{O}_S(2n,1)$ is very
ample so for $k$ big enough, $H^i(S, \omega_{f}^{\otimes k})=0$ $\forall\, i >0$.
Once again Leray spectral sequences give $H^i(S, \omega_S^{\otimes k})=H^i(\p,
\cal{R}_k)$ $\forall\, i, \; \forall\, k>1$.

Thus for $k$ big enough we get $\chi(\cal{R}_k)= h^0(S, \cal{O}_S(2nk,k))$. We can calculate
the latter tensoring the resolution (\ref{risoluzione}) with $\cal{O}(2nk,k)$, obtaining a new exact
sequence, which yields
\[\chi(\cal{R}_k)=\frac12 (-8 + 16k + 20n -    41kn + 41k^2n)\]
at least for big $k$.
We know that \[{k\choose2} K_f^2=\chi({R}_k)-\chi_f-4(2k-1)=41n{k\choose2}\] (see equation (\ref{chi})), then $K_f^2=41n$. Since $\operatorname{length}(\cal{F})=2n$, this is exactly the thesis of Theorem \ref{enunciato}.

\begin{rmk}When $n=1$, the surface $S$ can be seen as a complete intersection in
$\mathbb{P}^4$:
\begin{equation*}
\begin{matrix}
S & \longrightarrow & \mathbb{P}^4\\
(t_0,t_1),(x_0, \ldots, x_4) & \mapsto & (x_0, \ldots, x_4)
\end{matrix}
\end{equation*}
Its image obviously lies in the surface \[\tilde{S}=\{x_0 q_3 - x_2 q_2+x_3q_1=0, x_1 q_3 - x_3
q_2+x_4q_1=0\}.\] The morphism $S \fr \tilde{S}$ is in fact an isomorphism, as it has an inverse morphism. This
can be shown by calculation for $q_1, \,q_2, \, q_3$ general enough.\end{rmk}

\begin{rmk} This example clarify the dependence of the Horikawa index both on the geometry of the fibre and on the embedding of the fibre itself. For any $n \in \mathbb{N}$, the fibre over $(1,0)$ is trigonal, thus has positive Horikawa index. But for a generic choice of the polynomials involved, we find a suitable fibration for any $n \in \mathbb{N}$, thus the trigonal fibre is exactly the same, whether the index is $n$, thus the precise value of the index depends only on the embedding of the same fibre in the different fibration. 
\end{rmk}

\begin{es}
The previous example can be adapted to a more general ambient space, namely a normal rational scroll
$\mathbb{F}(a_0,a_1,a_2,a_3,a_4)$. In particular we can easily find examples of regular surfaces with any odd $p_g \geq 5$ such that $K_f^2=4\chi_f+1=4p_g+21$.\end{es}

Let us fix $\mathbb{F}=\mathbb{F}(a,a,0,0,0)$ for any $a \geq 0$. We can consider the Pfaffian
equations of a matrix very similar to the one of the previous example:

\begin{equation}
M=\left( \begin{matrix}
  & t_1 & t_0^a x_0 & x_2 & x_3\\
 & & t_0^{a+1} x_1 & t_0x_3 & t_0x_4\\
 & & & q_1 & q_2\\
 & & & & q_3
\end{matrix} \right)
\end{equation}
with Pfaffian equations

\begin{equation} \left \{\begin{array}{l}

c_1=t_0(t_0^a x_1q_3-x_3q_2+x_4q_1)\\
c_2=t_0^a x_0q_3-x_2q_2+x_3q_1\\
p_1=t_1q_3+t_0( - x_2x_4+x_3^2)\\
p_2=t_1q_2+t_0^{a+1}(x_0x_4-x_1x_3)\\
p_3=t_1q_1+t_0^{a+1}(x_0x_3-x_1x_2)

\end{array}\right.\end{equation}
where the $q_i$ have bidegree $(a,2)$. Like in Example \ref{ex}, we can find $q_i$ general enough to yield a
smooth surface with one single trigonal fibre.

We can show again that the surface is regular and that $\omega_S\cong\cal{O}_S(0,1)$ and
$H^0(\mathbb{F}, \cal{O}_{\mathbb{F}}(0,1)) \stackrel{\sim}{\fr}H^0(S,\omega_S)$. So $p_g=2a+5$.

\begin{es}We can modify the latter example in order to obtain even $p_g\geq 6$. If $p_g$ is even and the surface is regular, the degree of $\cal{R}_1$ is odd, so we look for an ambient space $\text{Proj(Sym}^n \cal{R}_1)=\mathbb{F}(a_0,a_1,a_2,a_3,a_4)$ with $\sum_i a_i$ odd. \end{es}

Let us fix $\mathbb{F}=\mathbb{F}(a,0,0,0,0)$ for any $a$ positive odd integer of the form $2d-1$. The matrix involved is

\begin{equation}
M=\left( \begin{matrix}

  & t_1^{d+1} & t_0^{2d} x_0 & x_2 & x_3\\
 & & t_0^{d+1} x_1 & t_0^dx_3 & t_0^dx_4\\
 & & & q_1 & q_2\\
 & & & & q_3
\end{matrix} \right)
\end{equation}
with Pfaffian equations

\begin{equation} \left \{\begin{array}{l}

c_1=t_0^d(t_0 x_1q_3-x_3q_2+x_4q_1)\\
c_2=t_0^{2d} x_0q_3-x_2q_2+x_3q_1\\
p_1=t_1^{d+1}q_3+t_0^d( - x_2x_4+x_3^2)\\
p_2=t_1^{d+1}q_2+t_0^{d+1}(t_0^{2d-1}x_0x_4-x_1x_3)\\
p_3=t_1^{d+1}q_1+t_0^{d+1}(t_0^{2d-1}x_0x_3-x_1x_2)

\end{array}\right.\end{equation}
where $q_1$ and $q_2$ have bidegree $(0,2)$ and $q_3$ has bidegree $(-1,2)$.

We can show again that the surface is regular and that $\omega_S\cong\cal{O}_S(0,1)$ and $H^0(\mathbb{F}, \cal{O}_{\mathbb{F}}(0,1)) \stackrel{\sim}{\fr}H^0(S,\omega_S)$. So $p_g=a+5=2d+4$.

\begin{rmk} The only missing values for $p_g$ are $0,\ldots, 4$. These can be obtained by a suitable modification of the above construction. \end{rmk}

\vskip1cm
\noindent Elisa Tenni\\
Dipartimento di Matematica ``F. Casorati''\\
Universit\`a degli Studi di Pavia\\
via Ferrata 1, 27100 Pavia, Italy\\
\textit{e-mail:} elisa.tenni\texttt{@unipv.it}
\end{document}